\documentclass[12pt]{article}
\usepackage{epsfig}
\usepackage[sectionbib]{natbib}
\newtheorem{theorem}{Theorem}
\newtheorem{lemma}{Lemma}
\newtheorem{definition}[theorem]{Definition}

\newcommand{\be}{\begin{eqnarray}}
\newcommand{\ee}{\end{eqnarray}}

\def\e{\hbox{E}}

\def\Min{\hbox{min}}
\def\Max{\hbox{max}}

\def\sign{\hbox{sign}}
\def\arg{\hbox{arg}}

\begin{document}
\title{\bf On two-sided p-values for non-symmetric distributions}
\author{Elena Kulinskaya\thanks{Statistical Advisory Service,
8 Princes Gardens, Imperial College, London SW7 1NA,
UK. Tel +4420 7594 3950.
e-mail: e.kulinskaya@ic.ac.uk}}
\date{October 5, 2008}\maketitle

\begin{abstract}
\noindent Two-sided statistical tests and p-values  are well defined only
when the test statistic in question has a symmetric distribution.
A new two-sided p-value  called  {\it conditional  p-value} $P_C$ is introduced here.
 It is closely related
to the doubled p-value and has an intuitive appeal. Its use is advocated
for both continuous and discrete distributions. An important  advantage of this p-value
is that equivalent 1-sided tests are transformed into  $P_C$-equivalent
  2-sided tests. It is compared to the widely used doubled and minimum likelihood p-values.
  Examples include the variance test, the binomial and the Fisher's exact test.
\\
\noindent {\bf keywords: }{two-sided tests, Fisher's exact test, variance test, binomial test, $F$ test,
minimum likelihood}
\end{abstract}\newpage

\section{Introduction}
Two-sided statistical tests are widely used and misused in numerous applications of statistics.
 In fact, some applied journals  do not accept papers quoting 1-sided
p-values anymore. Examples include The New England Journal of Medicine,
Journal of the National Cancer Institute
and Journal of Clinical Oncology
among others.

 Unfortunately, two-sided statistical tests and p-values  are well defined only
when the test statistic in question has a symmetric distribution. The difficulties with two-sided p-values
arise in a general case of a non-symmetric distribution, though they are more often commented on for
discrete distributions.

The most famous example is an ongoing discussion about how 2-sided p-values
 should be constructed for the Fisher's exact test .
 This discussion was started in 1935 by \citet{Fisher} and  \citet{Irwin}.
 Numerous developments of the next 50 years are summarised in \citet{Yates} and
discussion thereof. The more recent contributions include
several proposals based on an a modified UMPU test
\citet{Lloyd}, \citet{Dunne}, \citet{meulepas}. See also \citet{Agresti-1977}, \citet{Dupont}, \citet{Davis},
\citet{Agresti-1992}.
 A (far from exhausting) list of 9 different proposals is given
in \citet{meulepas}.  The problem is still not resolved.\\

Fisher advocated doubling the 1-sided p-value in his
letter to Finney in 1946 \citep [p.444]{Yates}. This doubled p-value is denoted by $P_F$.
Fisher's motivation was an equal prior weight of departure in either direction. Other arguments for doubling include
invariance under transformation of the distribution to a normal scale, and ease of approximation by the chi-square distribution
\citep{Yates}. One of the evident drawbacks of the doubling rule is that it may result in a p-value greater than 1. The doubled p-value is used in the majority of statistical software in the case of continuously distributed statistics and often in the discrete case. 
\\

The primary contribution of this article is the introduction of a new method of
defining two-sided p-values to be called \lq {\it conditional two-sided p-values}' denoted by $P_C$.
 The conditional p-value is closely related
to the doubled p-value and has an intuitive appeal. It is demonstrated that this new two-sided p-value has properties which make it a
definite improvement on currently used two-sided p-values for both
discrete and continuous non-symmetric distributions.\\

Another popular two-sided p-value for non-symmetric discrete distributions implemented in computer packages, R \citep{rrr} in particular, is adding the probabilities
of the points less probable than the observed (at both tails). This p-value is denoted by $P_{prob}$. This method
was introduced in \citet{F-H-1951}, and is based on \citet{NP-1931} idea of ordering multinomial
probabilities; this is called \lq the principle of minimum likelihood' by \citet{GP-1975}, see also \citet{george1990}.
\citet{HP-1965} were the first to use this p-value for Fisher's exact test.
Many statisticians objected to this principle. \citet{GP-1975}
 commented that \lq The minimum likelihood method can also lead to
absurdities, especially when the distribution is U-shaped, J-shaped, or simply not unimodal.' \citet{RA-1975} pointed out \lq This procedure is justified only if events of lower probability
are necessarily more discrepant from the null hypothesis. Unfortunately, this is frequently not true.'\\

The following example clearly demonstrates another unfortunate feature of this p-value. When a value of density is associated with a high 1-sided p-value at one tail,
the value at the opposite tail cannot be rejected even though it may have a very low  1-sided p-value.
\\

\begin{figure}[t]\label{fig:chi1}
  \vspace{-1cm}\caption{Two-sided variance test  with the statistic $X\sim \chi^2(5)$ .
     On the left plot, the density of $\chi^2(5)$
  distribution, with dotted/dashed lines illustrating the calculation of the $P_{prob}$ for $X=1$ and $X=0.5$.
  On the right plot, the power of the 5\%-level variance tests  based on the p-values
   $P_{prob}(x)$ (solid line),
   $P_F(x)$ (dashed line), $P_{C}^E(x)$ (dotted line), and the UMPU test (long-dashed line). The horizontal line at $0.05$ corresponds to the significance level. }
   \centerline{\psfig{figure=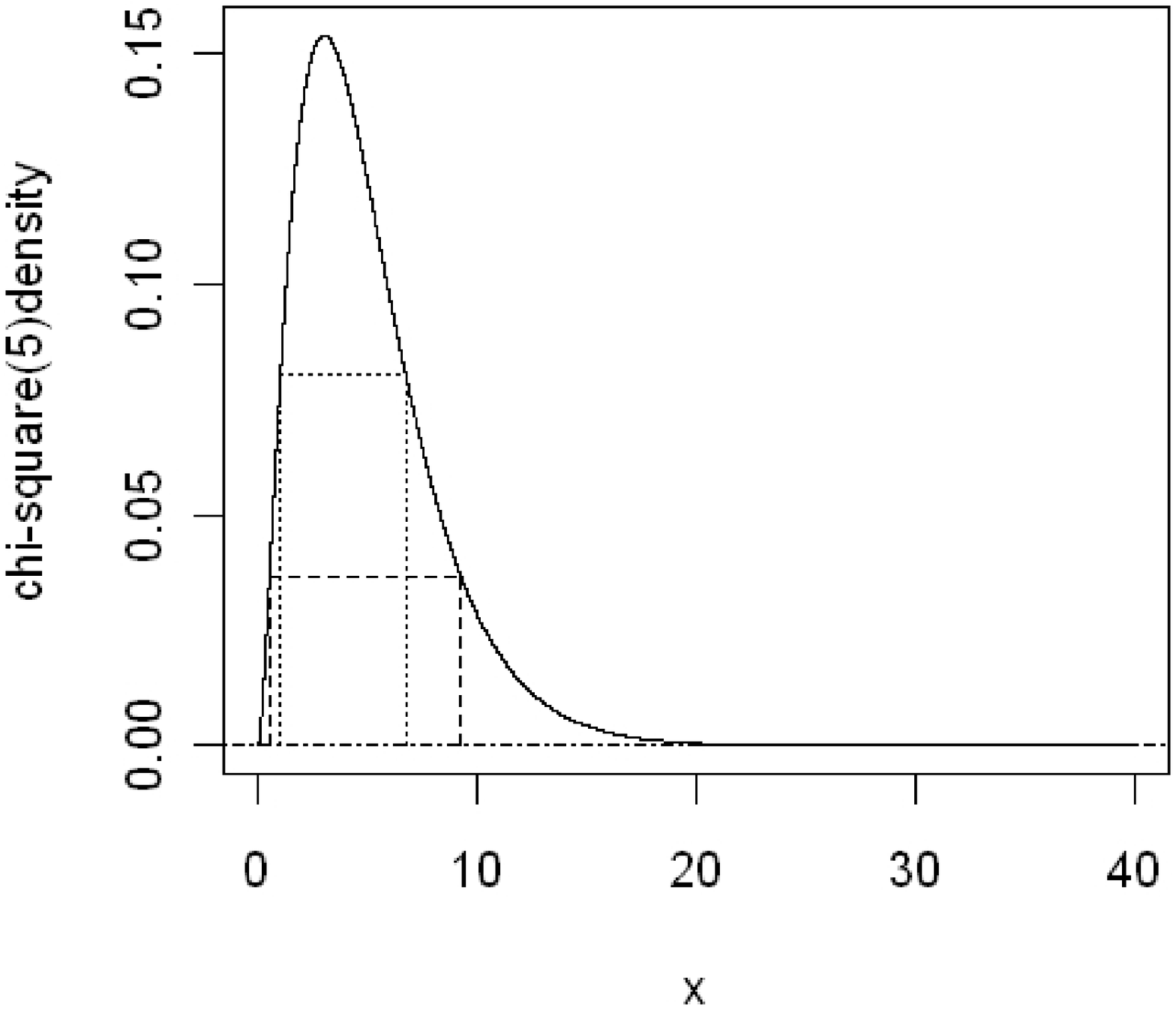,height=8truecm,width=8truecm}
    \psfig{figure=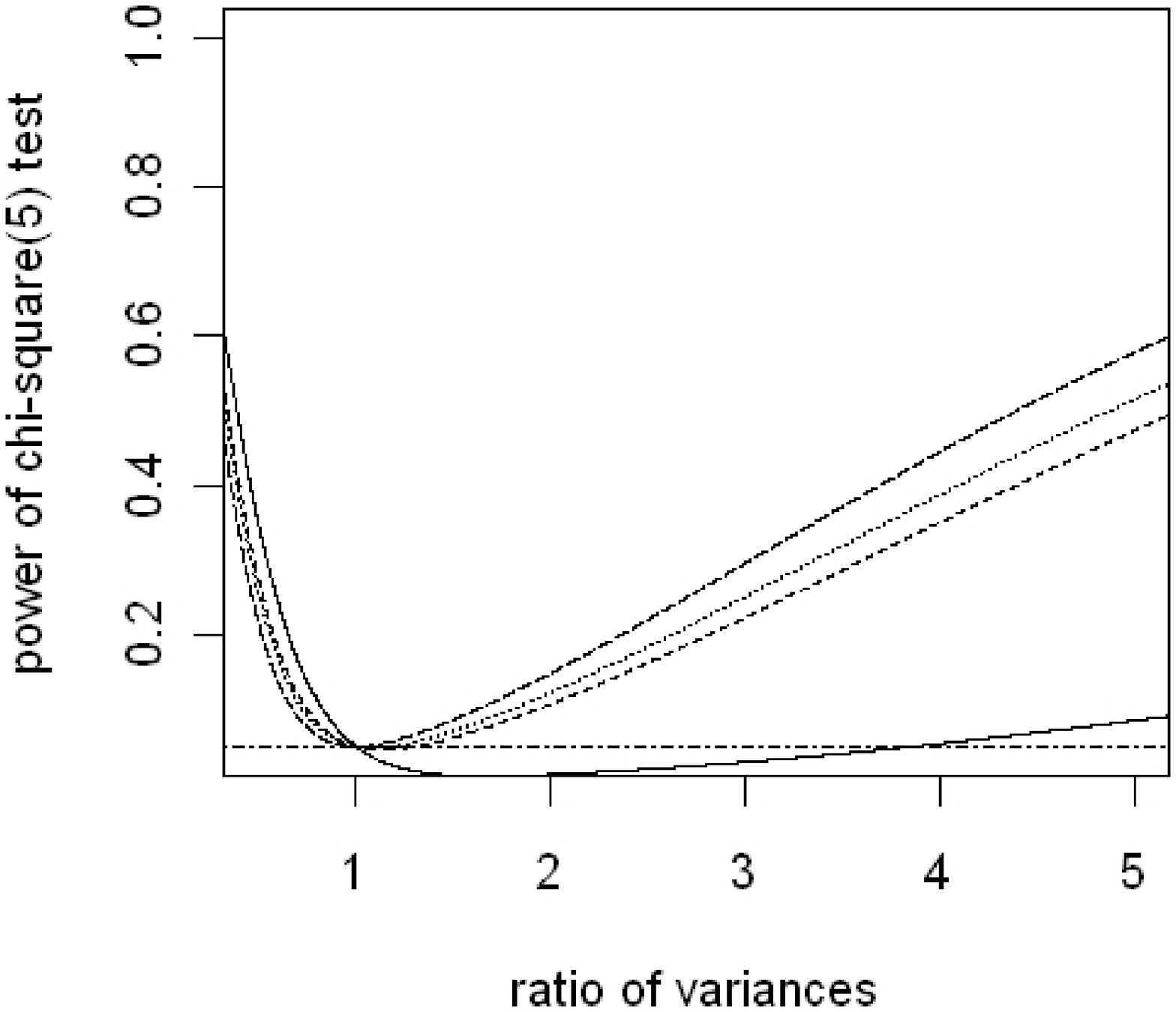,height=8truecm,width=8truecm}}
  \vspace{-.5cm}
\end{figure}

{\bf Example: Two-sided variance test based on the Chi-square distribution}
 Suppose we have 6 observations from a perfectly normal population and wish to test the null hypothesis that
the variance $\sigma^2=\sigma^2_0$ against a two-sided alternative.
The test statistic $X=(n-1)S^2/\sigma^2_0\sim (\sigma^2/\sigma_0^2)\chi^2(5)$, where $S^2$ is the sample variance. For $X=1$ (or $S^2=0.2$) the 1-sided
p-value on the left tail is 0.0374, the density is 0.0807,
the symmetric value on the right tail is $x'=6.711$, the 1-sided p-value  is 0.2431, see dotted lines on the
left plot of Figure \ref{fig:chi1}; similarly for $X=0.5$ ($S^2=0.1$) the density is 0.0366, the p-value on the left tail
is 0.0079, the symmetric value is 9.255, p-value=0.0993 (dashed lines on the same plot). It is very difficult to reject the null for small observed values.

Given critical values on the left and right tail $c_{L,\alpha}$ and $c_{R,\alpha}$, such that $\chi^2_{5}(c_{L,\alpha})+1-\chi^2_{5}(c_{R,\alpha})=\alpha$, the power of a two-sided variance test at level $\alpha$ is calculated as $\chi^2_{5}(\rho c_{L,\alpha})+1-\chi^2_{5}(\rho c_{R,\alpha})$, where $\rho=\sigma_0^2/\sigma^2$.
The power of four  $0.05$-level tests is plotted at the right plot of Figure~\ref{fig:chi1}. The test based on the $P_{prob}$ is evidently biased,
with
very low power for $\rho<1$, i.e. when $\sigma<\sigma_0$. The minimum value of the power is $0.01$.\\

The uniformly most powerfull unbiased (UMPU) test for this example has the critical region defined by $c_{L,.05}=0.989$ and $c_{R,.05}=14.37$ corresponding to critical levels $\alpha_L=.037$ and $\alpha_R=.013$ at the left and right tail, respectively. Finally, the generalized likelihood ratio (GLR) test is based on the statistic $\Lambda~=~[(X/n)\exp(1-(X/n))]^{n/2}$, and it is biased \citep [Example 23.5, p.882]{Kendall}. This is not exceptional; \citet{BarLev2002} showed that for  a continuous exponential family $F$ on the real line, the GLR and UMPU tests coincide if and only if, up to an affine
transformation, $F$ is either a normal, inverse Gaussian or gamma family.

 The new conditional 2-sided p-value $P_C$ is formally defined in the next section.
  The power of the tests based on the doubled and conditional p-value
for the chi-square example is also plotted in
Figure~\ref{fig:chi1}. They are much less biased, with minimum power of $0.045$ and $0.048$ respectively.
\\

The formal definition of the conditional 2-sided p-value $P_C$ and the comparison of its
properties  to those of the doubled p-value and the $P_{prob}$ for a case of continuous distributions is given
in Section 2,
and for discrete distributions (binomial and hypergeometric) in Section 3. Discussion is in Section 4.
The use of the conditional 2-sided p-value $P_C$ is advocated
for both continuous and discrete distributions. An important  advantage of this p-value
is that equivalent 1-sided tests are transformed into  $P_C$-equivalent
  2-sided tests.

\section{Two-sided p-values for continuous asymmetric distributions}
Consider a general case of a statistic $X$ with a strictly increasing continuous
 null distribution $F(x)$  with continuous density $f(x)$. For an observed value $x$ of $X$,  one-sided p-value on the left tail is defined as $P(X'\leq x|X=x)=F(x)$, where $X'\sim F(x)$ independent from $X$. Similarly, on the right tail the p-value is
 $P(X'\geq x|X=x)=1-F(x)$. Denote by $A$ a generic location parameter chosen to  separate the two tails of the distribution $F$. Particular examples include the mean $E=\e(X)$, the mode $M=\arg\; \sup_x f(x)$,
   or the median $m=F^{-1}(1/2)$. 
  Which parameter should be used to separate the two tails depends on the context;
  the mean  seems to be the most appropriate when a test statistic is based on an estimate of a natural parameter in an
  exponential family, as is the case with binomial or Fisher's exact test.
  General theory below  is applicable regardless of the parameter chosen, though the details of examples may differ.
  Interestingly, it does not matter much
  for the most important non-symmetric discrete distributions: the mean when attainable coincides with the mode (or one of the two modes) for Poisson, binomial and hypergeometric distributions.  The latter two distributions are discussed in Section 3.\\

 \begin{definition}\label{def:pweighted}
 Weighted two-tailed p-value centered at $A$ with  weights $w=(w_L,w_R)$ satisfying $w_L+w_R=1$ is defined as
 \begin{equation}\label{eq:weighted}P^A_{w}(x)=\Min \bigl (\frac{F(x)}{w_L}|_{x< A},\frac{1-F(x))}{w_R}|_{x>A},1\bigr ).\end{equation}
  \end{definition}

 Doubled p-value denoted by $P_F^A(x)$  has weights $1/2$. Without loss of generality assume that $A>m$.  Then the doubled p-value $P_F^A(x)$ is equal to $2F(x)$ for $x<m$,
1 for $m\leq x\leq A$, and $2(1-F(x))$ for $x>A$.  Thus the doubled p-value is not continuous at $A$ unless $m=A$, its derivative is also discontinuous at $m$. \\

Similarly, a weighted p-value $P_w^A(x)$ is continuous at $A$
iff $w_L/w_R=F(A)/(1-F(A))$ and an additional requirement of $P_{w}^A(A)=1$ results in
 $w_L=F(A)$ and $w_R=(1-F(A))$ arriving at the next definition.


 \begin{definition}\label{def:conditional}
 Conditional 2-sided p-value centered at $A$ is  defined as \begin{equation}\label{eq:pcond}
 \begin{array}{ll}P_C^A(x)&=P^A_{\{F(A),1-F(A)\}}(x)\\&=P(X'\leq x|X=x\leq A)+P(X'\geq x|X=x\geq A).\end{array} \end{equation}
\end{definition}
This is a smooth function of $x$ (but at $A$),
with a maximum of $1$ at $A$. It strictly increases for $x<A$ and decreases for $x>A$. The conditional p-value is conceptually close to the doubled p-value, the only difference being that the two tails are
 weighted inversely
proportionate to their probabilities. This results in  inflated p-values on the thin tail, and
deflated p-values on the thick tail when compared to the doubled p-value. When the tails are defined in respect to
 the median, the two p-values coincide: $P_F^m(x)=P_C^m(x)$. Thus conditional p-value is equal to the usual doubled
 p-value for a symmetric distribution.
It is easy to see that under the null hypothesis $P_C^A(x)$ is uniformly distributed on $[0,1]$ given a particular tail, i.e.
$P_0(P_C^A(X)\leq p|X\leq A)=p$,
similar to a 1-sided p-value.   \\


There is an evident connection between a choice of a two-sided p-value and a critical
region (CR) for a  two-sided test at level $\alpha$. A CR is defined through critical values corresponding to probabilities
$\alpha_1=w_L\alpha$ and $\alpha_2=w_R\alpha$, with the weights of the two tails
$w_L+w_R=1$. It can  equivalently be defined through a weighted p-value as $\{x: P_w^A(x)<\alpha\}$. The doubled p-value corresponds to $w_L=w_R=1/2$. The conditional p-value is equivalent to the choice  $w_L=F(A)$, $w_R=1-F(A)$.\\

For a two-sided test, critical values  $c_{L,\alpha}$ and $c_{R,\alpha}$ satisfy $F(c_{L,\alpha})=w_L\alpha$ and $1-F(c_{R,\alpha})=w_R\alpha$. Thus $w_L=w_L(\alpha)=F(c_{L,\alpha})/\alpha$. Define  $A=A(\alpha)=F^{-1}(F(c_{L,\alpha})/\alpha)$. Then the CR is $\{x: P^A_C(x)<\alpha\}$. 
 Therefore any 2-sided test,  a UMPU test inclusive, is a test based on conditional p-value centered at some $A=A(\alpha)$.  Conversely, if the $A$ value is chosen to be independent of $\alpha$, the resulting test is, in general, biased. Since an independence from level $\alpha$ is  a natural requirement for a p-value, some bias cannot be escaped. \\
 \begin{lemma}
For a one-parameter exponential family $F(x,\theta)$, a  two-sided level-$\alpha$ test based on the conditional p-value $P_C(A)$ is less biased in the neighborhood of the null value $\theta_0$ than the standard equal tails test based on the doubled p-value whenever $F(A)\in (1/2,w^*_{L,\alpha}]$, where $w^*_{L,\alpha}$ is the weight at the left tail of the UMPU test.
\end{lemma}
{\bf Proof} Denote  test critical function by $\phi(x)$. This is an indicator function of the CR, so $\e_0[\phi(X)]=\alpha$, and the power is $\beta(\theta)=\e_{\theta}[\phi(X)]$. Without loss of generality $X$ is the sufficient statistic. The derivative of the power is \citep[p. 127] { Lehmann}
\begin{equation}\label{eq:deriv}\beta'(\theta)=\e_\theta[X\phi(X)]-\e_\theta(X)\e_\theta[\phi(X)]\end{equation}
For an UMPU test $\beta'(\theta_0)=0$. For a test with weight $w_L$ at the left tail,
$$\beta'(\theta_0)=\int_{-\infty}^{F^{-1}(\alpha w_L)}xdF+\int^\infty_{F^{-1}(1-\alpha(1- w_L))}xdF-\alpha E.$$
For $\alpha<1$, this is strictly decreasing function of $w_L$ equal zero at $w^*_{L,\alpha}$. When $1/2< w^*_{L.\alpha}$, any $w_L\in (1/2,w^*_{L,\alpha}]$ provides positive values of $\beta'(\theta_0)$, and when $1/2> w^*_{L.\alpha}$, the values of $\beta'(\theta_0)$ are negative; in any case the gradient is the steepest  and the bias is the largest at $1/2$, as required.\\

\begin{figure}[t]
  \vspace{-1cm}\caption{Bias of the $P_F(x)$-based variance test at 5\% level  (dashed line), and $P_{C}^E(x)$-based test (dotted line)
   in the 1-sample case ($\chi^2$-test, left plot) and in the 2-sample case with $n_1=6$ (F-test, right plot).  }
   \label{fig:bias}
    \centerline{\psfig{figure=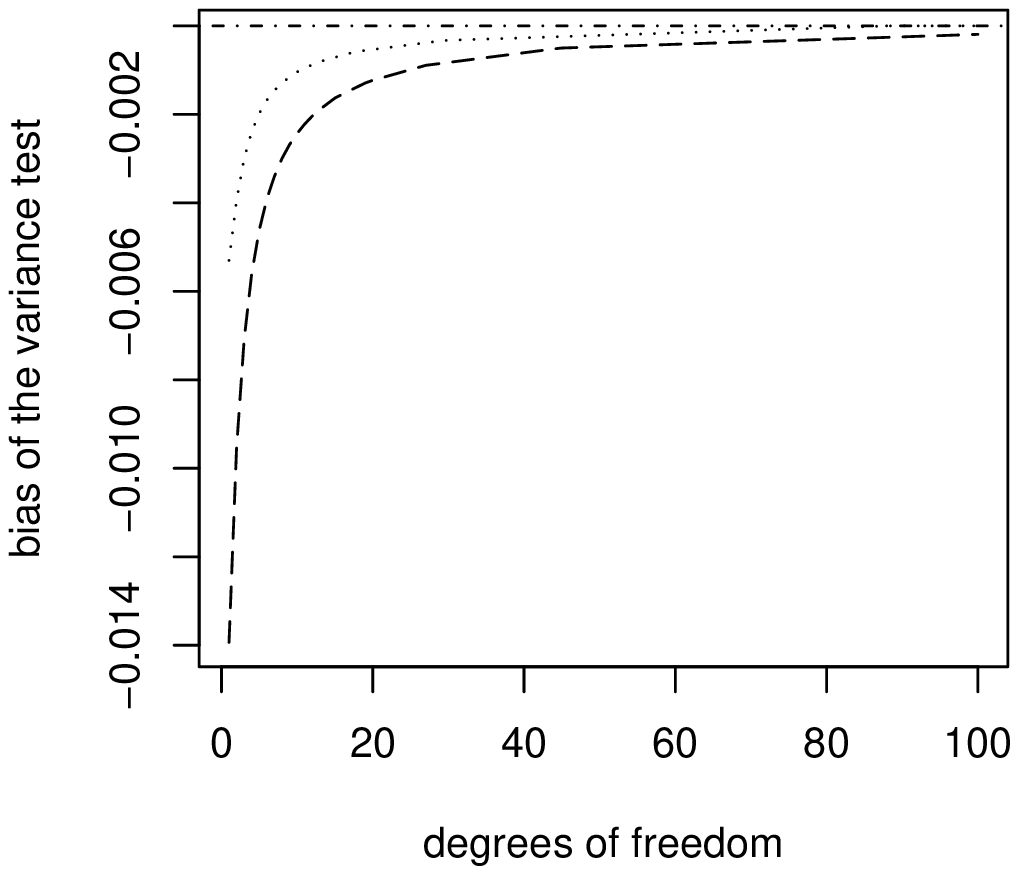,height=8truecm,width=8truecm}
    \psfig{figure=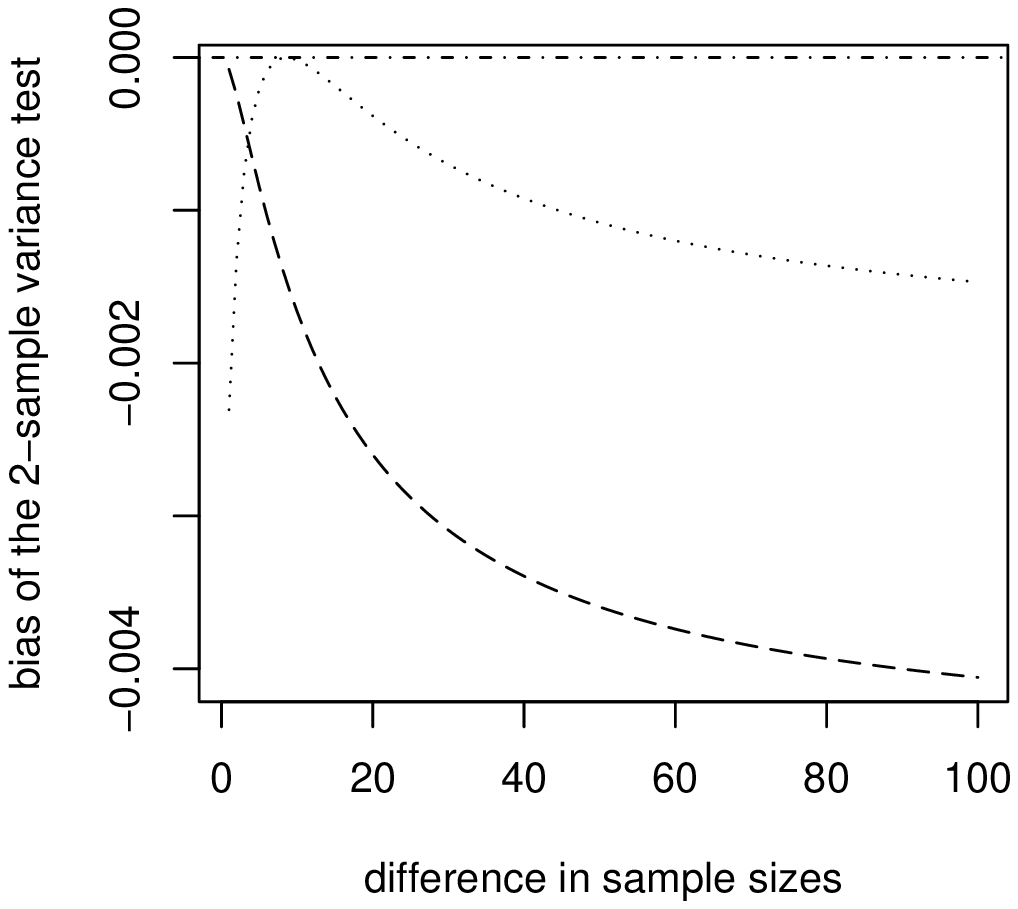,height=8truecm,width=8truecm}}
      \vspace{-.5cm}
\end{figure}

Lemma 1 provides a sufficient condition for the $P_C^E$-based test to be less biased than the equal tails test, but this condition is not necessary. This condition holds for the $\chi^2$ distribution, and the variance test
 based on $P_C^E(x)$  is uniformly (in $n$) less biased then the test based on the doubled p-value, left plot of Figure~\ref{fig:bias}.
The doubled p-value based test is asymptotically UMPU, \citet{Shao}, and so is the $P_C$-based test.  In the two-sample case, the equal-tails $F$-test of the equality of variances is UMPU for equal sample sizes, and the $P_C$-based test is  less biased when the ratio of sample sizes is larger than 1.7 (starting from $n=6$), whereas lemma 1 holds for even more unbalanced sample sizes with the ratio of 2.5 or above, right plot of Figure~\ref{fig:bias}.\\

Finally, consider the minimum likelihood p-value.
\begin{definition}\label{def:p_prob} Minimum likelihood p-value is $P_{prob}(x)=P(f(X)\leq f(x))$. \end{definition}
$P_{prob}(x)$ reaches 1 at the mode, and $P_{prob}(A)<1$ whenever $A\not = M$. It is not a unimodal function of $x$ when the density
$f(x)$ is not unimodal. In a case of a unimodal distribution, for a pair of conjugate points $(x,x'):x<M<x',f(x')=f(x)$, it is calculated as
$P_{prob}(x)=F(x)+1-F(x')$.  It has a $Unif(0,1)$ distribution under the null.

When used to define a test, the acceptance region defined as $\{x:P_{prob}~(x)>~\alpha\}$ contains the points
with the highest density, and is therefore of minimum length.
Inverting this test results in the shortest
confidence intervals, see \citet{Sterne} for the binomial and \citet{Baptista} for the hypergeometric distribution.
It is also related to Bayes shortest posterior confidence intervals,
see \citet{WT-1971} for the intervals for the standard deviation $\sigma$ and
the ratio of variances in normal populations, based on
inverse chi and $F$ distribution, respectively.\\

The next three examples clarify the properties of the  conditional p-value $P_C$ in comparison to $P_{prob}$.\\

{\bf Example: Triangular distribution}

Suppose that the null density is given by $f(x)=2(x+a)/[a(a+b)]$ for $-a\leq x\leq 0$ and $f(x)=2(b-x)/[b(a+b)]$ for $0\leq x\leq b$. The mode $M=0$, and $F(0)=a/(a+b)$. Then for $x<0$, $P_{prob}=F(x)/F(0)$ and for $x>0$, $P_{prob}=(1-F(x))/(1-F(0))$, \citet{george1990}. Thus, $P_C^M(x)=P_{prob}(x)$. This is the only unimodal distribution for which this equality holds as it requires the linearity of the density $f(x)$.

 {\bf Example: Uniform distribution}

Consider a $Unif(0,1)$ distribution. This is a symmetric distribution with $E=m=1/2$, and $P_C=P_F=2x$ for $x\leq 1/2$,
 and $P_C=P_F=2(1-x)$
for $x\geq 1/2$, whereas  $P_{prob}\equiv 1$ for all values of $x\in[0,1]$. \\
This example shows the cardinal difference between
the two p-values. $P_C$ acknowledges unusual values of $x$ at the ends of the interval, and the $P_{prob}$ does not.
This is a somewhat extreme example, because the uniform distribution has a whole interval of modes. The next example
deals with a unimodal distribution but shows exactly the same properties of the respective p-values. \\

{\bf Example: Left-truncated normal distribution.}

\begin{figure}[t]
  \vspace{-1cm}\caption{Plot of $P_{prob}(x)$ (solid line), $P_F^E(x)$ (dashed line), and $P_{C}^E(x)$ (dotted line)
   for the $\chi^2(5)$ distribution (left plot) and for a standard normal distribution truncated
  at $-0.5$ (right plot). The plotted doubled p-value $P_F^E(x)$ is not truncated at 1. }
  \label{fig:trunc}
    \centerline{\psfig{figure=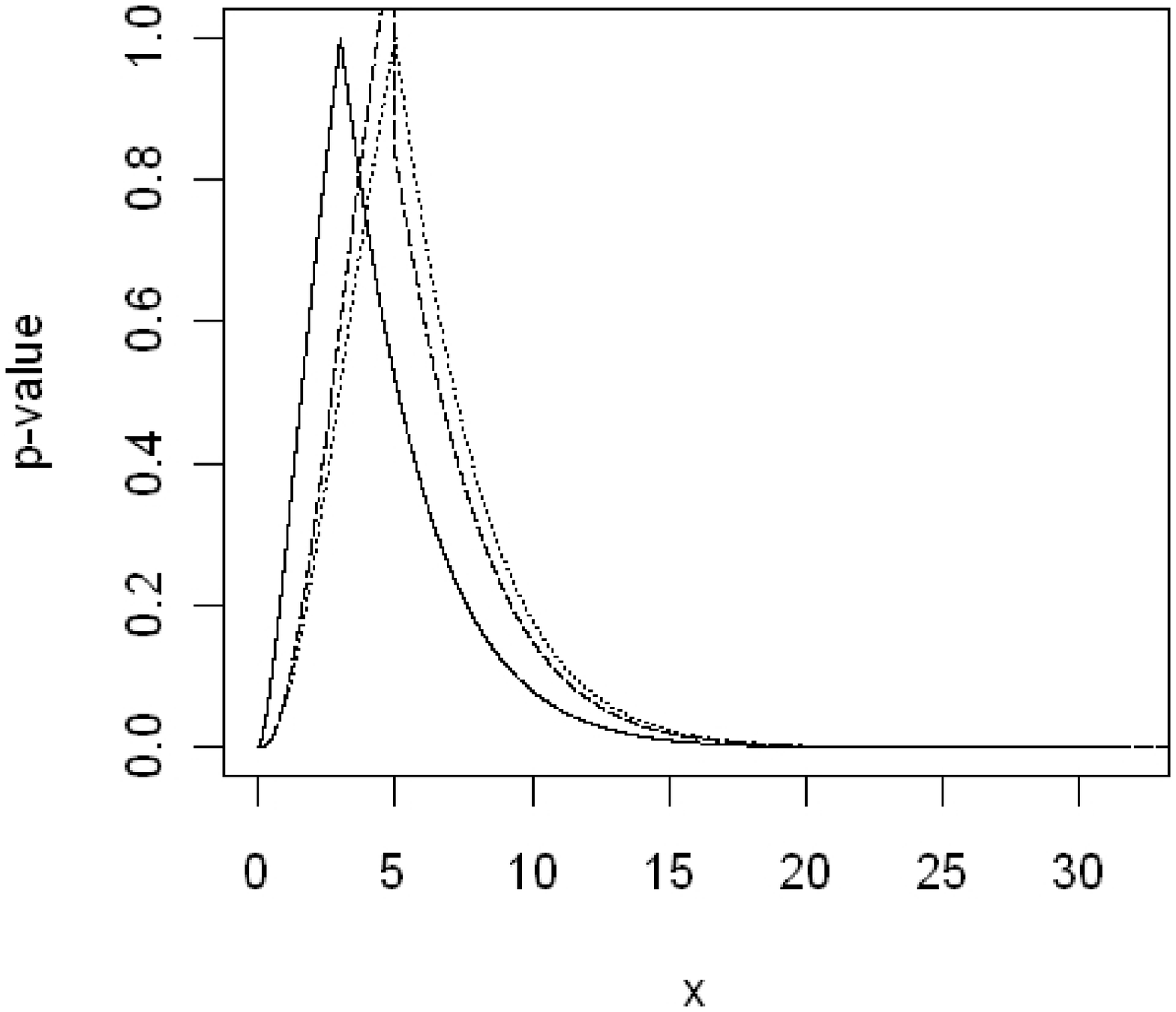,height=8truecm,width=8truecm}
    \psfig{figure=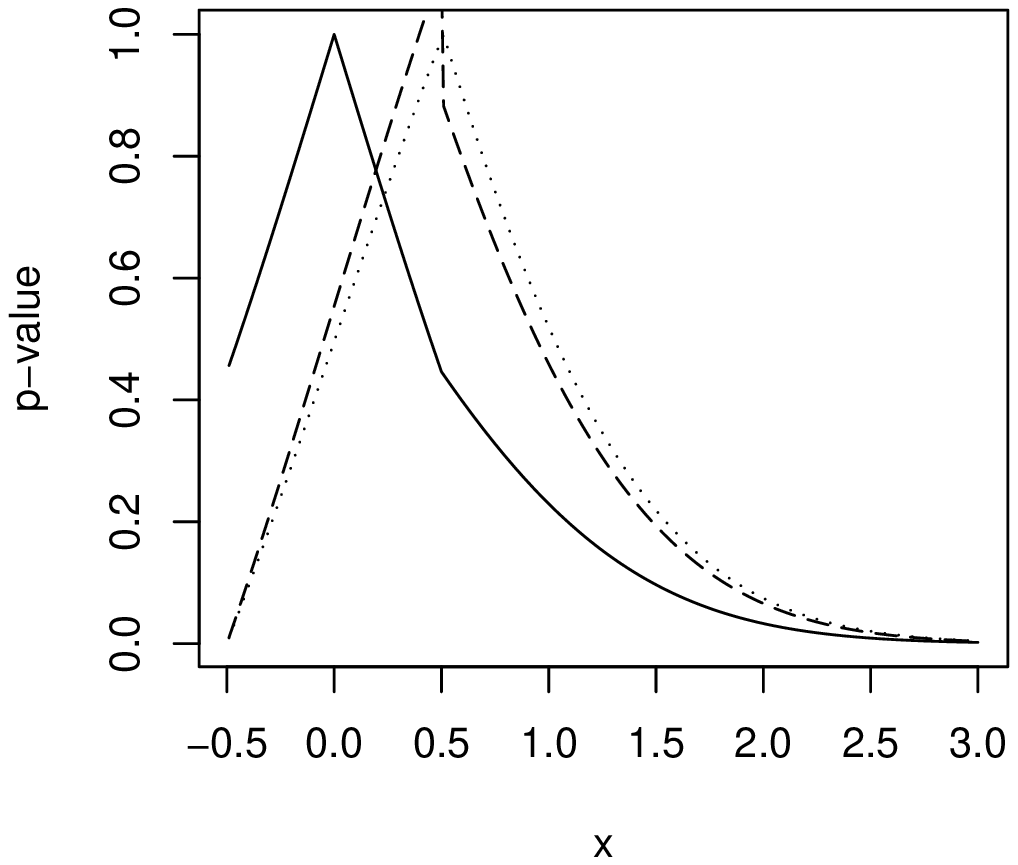,height=8truecm,width=8truecm}}
  \vspace{-.5cm}
\end{figure}

Denote the standard normal distribution function and density by $\Phi(x)$ and $\phi(x)$, respectively.
Consider a left-truncated at $-L<0$ normal distribution $G_L(x)=(\Phi(x)-\Phi(-L))/(1-\Phi(-L))$ defined for  $x\geq -L$.
The mode is at zero. Then $P_{prob}(x)=2G_L(-|x|)+1-G_L(L)$ for $-L\leq x\leq L$, and $P_{prob}(x)=1-G_L(x)$ for $x>L$.
$P_{prob}$ reaches 1 at $0$,
and $P_{prob}(\pm L)=1-G_L(L)$ is continuous at $L$, but its derivative is not continuous at $L$.
The mean is $E=E(L)=\phi(-L)/(1-\Phi(-L))$, and the conditional p-value $P_C(x)$ reaches 1 at $E(L)$.
 An example for $L=0.5$ is plotted in the right plot in Figure \ref{fig:trunc}.
For this example $E(L)=0.509$ and the weight of the left tail is $w_L=0.558$. The main difference between the two p-values is that $P_{prob}\geq 1-G(L)$ at the left tail, so even the low values of $x$
in the vicinity of  $-L$ have rather high p-values. On the other hand, $P_C$ is very close to zero for these values,
recognizing that it is rather unusual to get close to $-L$. It seems that a  small  two-sided p-value at the left tail
makes more sense. \\

The above two examples show the properties of the $P_C$ which are perhaps clear from its definition: it compares a value $x$
to other values at the same tail. On the other hand, $P_{prob}$ depends on the values at both tails.
The same circumstances arise in the variance test example which was introduced in the Introduction.\\

{\bf Example: variance test based on the Chi-square distribution (continued)}
Recall, that for $X=0.5$ ($s_n^2=0.1$)  the 1-sided p-value is .0079, and the value of $X$ with equal density
is $X'_p=9.256$ with the 1-sided p-value of .0993. The mean $E=5$, and
the conditional p-values are $P_C^E(0.5)=0.0135$, and $P_C^E(9.256)=0.239$, the weight of the left tail is $w_L=F(E)=0.584$.
The value with the same conditional p-value on the opposite tail is
$X'_C=F^{-1}(1-(1-w_L)P_{C}(0.5))=16.48$ with the 1-sided p-value of 0.0056. Clearly, $X'_C$ is more comparable to $X$
than the value $X'_p$. The three p-values are plotted at the left plot in Figure~\ref{fig:trunc}.

The power of the three tests and of the UMPU test (all at 5\% level) is shown in the right plot in Figure~\ref{fig:chi1}. The UMPU test is the  conditional test with $A=6.403$, corresponding to the weight $w^*_L=0.731$. All three tests are biased,  with the bias $B$ defined as the minimum difference between the power and level
being $B_F=-0.0046$ for the doubled and $B_C= -0.0020$ for the conditional test. This agrees with Lemma 1.
The doubled test is slightly less powerful on the right, and slightly more on the left.
 The test based on $P_{prob}$ has very large bias and such low power for the alternatives $\sigma<\sigma_0$, that it
does not deserve to be called a two-sided test. \\

The main difficulty associated with the two-sided tests is that two equivalent 1-sided tests may result
in distinct 2-sided tests. For the variance test example
  the tests based on $|s^2_n-\sigma_0^2|$ and on $|\log(s^2_n/\sigma_0^2)|$ are not equivalent.
  Let $D(x,A)$ be a measure of distance from $A$. It imposes
an equivalence of points at two sides of $A$:  each value $x<A$ has an equidistant value $x'_D: D(x,A)=D(x',A)$.
Two-sided tests based on $|X-A|$ and $D(X,A)$
 are not equivalent, generally speaking,  because for $x<A$, the equidistant value  $x'_D\not = 2A-x$.
 This results in different rejection regions for the two tests.
 The main advantages of the conditional p-value $P_C^A(x)$
 are given in the next Lemma.

 \begin{lemma}
  (i) For  a strictly increasing function $T(x)$, the conditional p-value
  $P_C(T(x)|T(A))=P_C^A(x)$.\\
 (ii) Suppose $D(x,A)$ strictly decreases for $x<A$ and strictly
increases for $x>A$, and $D(A,A)=0$. Define the conditional p-value for the distance $D(x,A)$ as
 $P_C(D(x,A))=P(D(x',A)\geq D(x,A)|X=x\leq A)+P(D(x',A)\geq D(x,A)|X=x\geq A)$.
 Then $P_C(D(x,A))=P_C(|x-A|)$.
 \end{lemma}
 The first statement of the lemma easily follows from the definition of $P_C^A(x)$, and for the second
 statement take $T(x)=D(x,A)\sign(x-A)$. This is a strictly increasing function of $x$, and the proof follows from part (i).

 The first part of the lemma ensures that equivalent 1-sided tests are transformed into  $P_C(x)$-equivalent
  2-sided tests. The second part states that the 2-sided tests based on any measure of distance from $A$ are $P_C(x)$-equivalent.
  This is true because  the conditional p-value ignores
 any equivalence between the points at different tails. \\

\section {Discrete distributions}

 In this section the 2-sided conditional p-value $P_C$ is defined for a discrete distribution.
 It is also compared to $P_F$ and $P_{prob}$ for two important cases: binomial and hypergeometric
 distributions. \\

The definition of the conditional p-value $P_C$ (\ref{eq:pcond}) is also applicable in a discrete case, but it
may require  a modification when the value $A$ is attainable.
Since the value $A$ belongs to both tails, the previously defined weights of the tails
 $w_L=P(x\leq A)$ and $w_R=P(x\geq A)$ add up to $1+P(A)>1$.
The modified weights of the tails are $w_L^{A(m)}=P(x\leq A)/(1+P(A))$ and $w_R^{A(m)}=P(x\geq A)/(1+P(A))$. This modification is akin to continuity correction.
The formal definition of $P_C(x)$ is
\begin{definition}Conditional two-sided p-value for a discrete distribution is
\begin{equation}\label{eq:pcondd}P_C^A(x)= \frac {P(X\leq x)}{w_L}|_{(x< A)}
+1|_{(x=A)}+\frac{P(X\geq x)}{w_R}|_{(x>A)}, \end{equation}
 where  the weights are $w_L=P(x\leq A)$ and $w_R=P(x\geq A)$.
 Modified conditional p-value $P^{A(m)}_C(x)$ is defined with weights  $w_L^m=P(x\leq A)/(1+P(A))$ and $w_R^m=P(x\geq A)/(1+P(A))$
  in equation (\ref{eq:pcondd}). \end{definition} Two definitions coincide when the value $A$ is not attainable.
In a discrete symmetric case when $A=E=m$ is an attainable value the values of $P^m_C(x)=P_F(x)$ are doubled 1-sided values,
and the values of $P_C(x)$ are $(1+P(A))$ times smaller, and the $P_C(x)$-based test is therefore more liberal.
The conditional p-value has a mode of 1 at $A$ when this
value is attainable, and two modes of 1 at the attainable values above and below $A$
when $A$ is not an attainable value. It has  discrete uniform distribution when restricted
to values at a particular tail, though not overall. In what follows we consider the case of $A=E$,
and use the notation $P_C(x)=P_C^E(x)$.
\subsection{Binomial distribution}
For $Binom(n,p)$ distribution  the mode is $M=\lfloor {(n+1)p}\rfloor=\lfloor{E+p}\rfloor$.
When $(n+1)p$ is an integer, $M=(n+1)p$ and $M-1$ are both modes,
and the mean $E=np\in(M-1,M)$ is unattainable. When  $E$ is an integer, $M=E$. In all cases the distance
$|M-E|<1$. The median is one of $\lfloor {np}\rfloor$ or $\lfloor {np}\rfloor\pm 1$. \\
Consider first the symmetric case  $p=0.5$.  For odd $n$, the value $(n+1)p$ is an integer, both tails of the distribution  have weight 0.5
and $P_C(x)=P_{prob}(x)$. For even $n$, the mean $E=np$ is an integer, $w_L>0.5$, but $w_L^m=0.5$.
Unmodified version $P_C(x)$ is symmetric
at $E$ with values $P_C(x)<P_{prob}(x)$ for $x\not = E$. The modified version $P^{(m)}_C(x)=P_{prob}(x)$.\\

Statistical packages differ in regards to the 2-sided p-values for the binomial test: R \citep{rrr} uses $P_{prob}(x)$ and StatXact (www.cytel.com) uses the doubled p-value.\\

\begin{figure}[t]
  \vspace{-1cm}\caption{Plot of $P_{prob}$ (solid line, circles), $P^m_C$ (long dash, filled circles),  $P_{C}$ (doted line, squares) and $P_F$ (dashed line, triangles) for
   $Binom(10,0.2)$
  distribution  (left), and  $Binom(11,0.2)$(right). On the right plot $P_C^m=P_C$.}
    \centerline{\psfig{figure=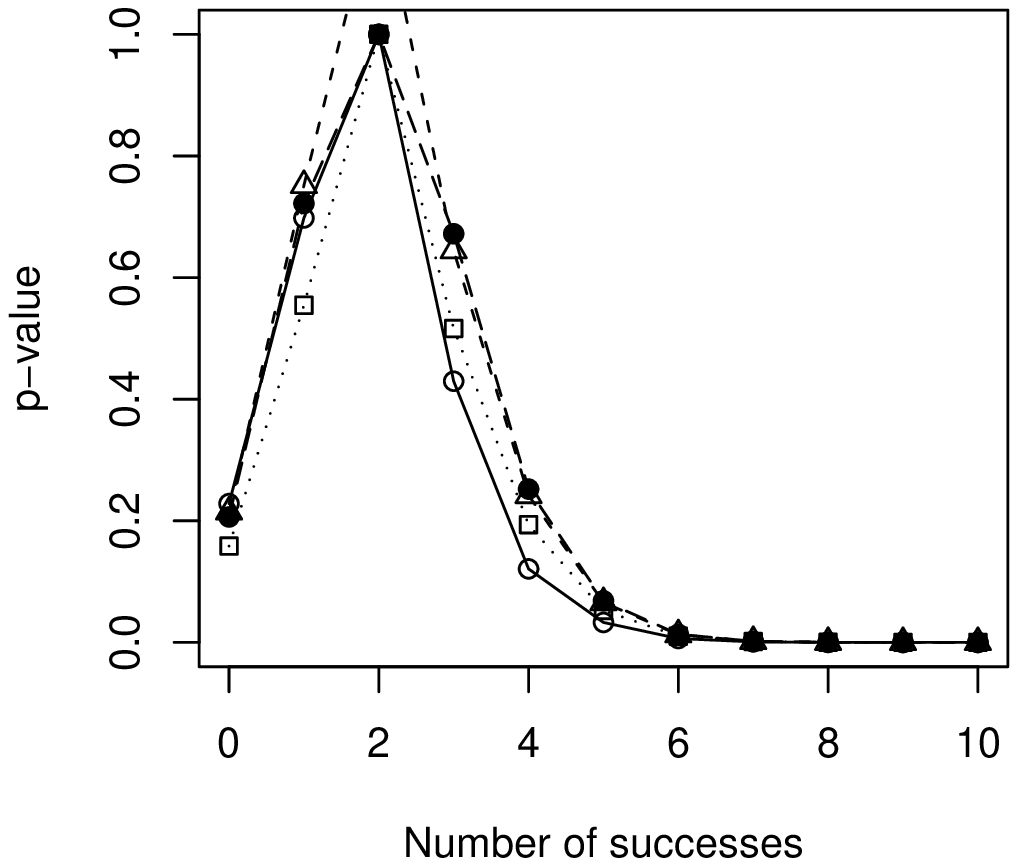,height=8truecm,width=8truecm}
    \psfig{figure=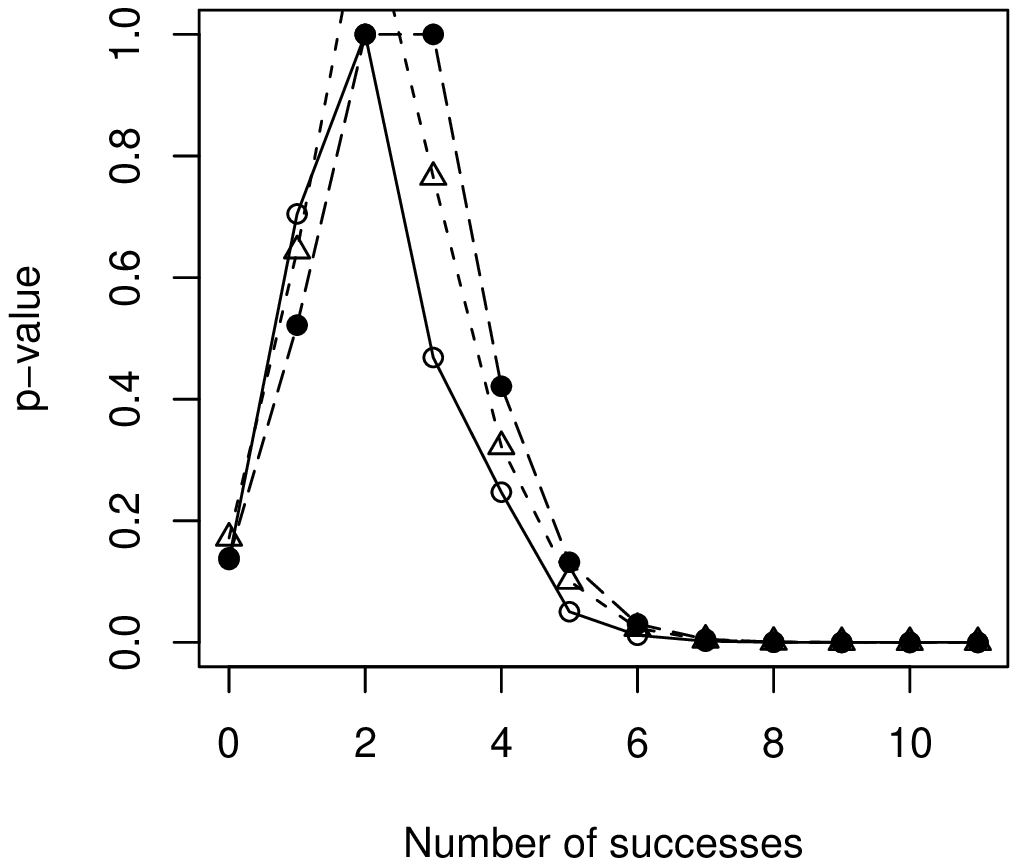,height=8truecm,width=8truecm}} \label{fig:d}
  \vspace{-.5cm}
\end{figure}

The three p-values, $P_C$, $P^{(m)}_{C}$, and $P_{prob}$  are plotted in Figure~\ref{fig:d}
 for $p=0.2$  and two values of $n$,  $n=10$ and  $n=11$.
 In the first case $E=2$ is an attainable value.
It can be seen that $P_C^{(m)}>P_C$ on the left plot. The weight of the left tail
 is $w_L=0.678$ vs $w_L^m=0.521$. Consequently, $P_C^m(x)=1.3P_C(x)$ for all $x$ but $E$.
 Modified conditional p-value $P^{(m)}_{C}$ is considerably
closer to $P_{prob}$ at the left tail,  and $P_{prob}<P_{C}<P^{(m)}_C$ at the right (thin) tail.
In fact, in this example for $n=10, p=0.2$, $P_{prob}$ provides exact 1-sided p-values for $x\geq 4$,
$P_{C}(x)=1.60P_{prob}(x)$ and $P^{(m)}_{C}(x)=2.09P_{prob}(x)$ for $x\geq 4$. So $P_{prob}(5)=0.033$, $P_F(5)=0.066$, $P_C(5)=0.052$, and
$P_C^m(5)=0.068$. The two-sided binomial test  as programmed in R uses $P_{prob}$ and would reject the null
hypothesis of $p=0.2$ at 5\% level given an observed value of 5, whereas a test based on the doubled or conditional p-value would not reject.
The same thing may happen for much larger values of $n$. For example, for $n=101$,
 $p=0.1$ and the observed value of $x=17$ the values are
$P_{prob}=0.030$ , $P_F=0.06$ and $P_C=P_C^m=0.052$.\\

For $n=11$ (right plot) $E=2.2$ is not attainable. $P_C=P_C^m$ has two modes at 2 and 3. Here $w_L=0.617$ and $P_C(x)=1.62PF(x)$ for $x\leq 2$, whereas $P_C=2.61(1-F(x-1))$ for $x\geq 3$.

Typically, $P_C(x)<P_{prob}(x)$ at the thick tail, and $P_C(x)>P_{prob}(x)$ at the thin tail.
Even for large $n$ the difference between $P_C$ and $P_{prob}$ is rather large. For example,
for $n=101$ and $p=0.1$ the values are $P_C(17)=0.052$ and $P_{prob}(17)=0.030$ in comparison to the 1-sided p-value of
$0.023$.

For the binomial distribution the weight of the tails converges to 0.5 rather slowly, and $P_C(x)\to P_F(x)$, see Table \ref{tab:binom}.
Whenever the mean is attainable,  the weight of the thin right tail is also more than 0.5,
  and $P_C(x) <P_F(x)$. If E is not attainable, the weight of the thin tail is less than 0.5,
  and then $P_C>P_F$. This is always true for $P^m_C(x)$. The distribution is more symmetric when the mean is
  attainable. Otherwise even for $n=1001$, the weight of the left tail is $w_L=0.522$ for $p=0.1$.

  \begin{center}
 \begin{table}[phtb]
\begin{small}
\begin{tabular}{|r||c|c|c|c|c|c|}
\hline
 & \multicolumn{3}{|c|}{$p=0.1$} & \multicolumn{3}{|c|}{$p=0.2$} \\
 \hline
 $n$ &   $w_L$   &   $w_L/w_R$    &   $w_L^m$  &   $w_L$   &   $w_L/w_R$    &   $w_L^m$   \\
\hline
10  &   0.736   &   1.130   &   0.531     &   0.678   &   1.086   &   0.521      \\
11  &   0.697   &   2.304   &   0.697    &   0.617   &   1.614   &   0.617     \\
20  &   0.677   &   1.113   &   0.527    &   0.630   &   1.070   &   0.517      \\
21 &   0.648   &   1.844   &   0.648    &   0.586   &   1.416   &   0.586     \\
50  &   0.616   &   1.083   &   0.520   &   0.584   &   1.049   &   0.512      \\
51 &   0.598   &   1.485   &   0.598    &   0.556   &   1.250   &   0.556      \\
100 &   0.583   &   1.063   &   0.515   &   0.559   &   1.036   &   0.509      \\
101 &   0.570   &   1.325   &   0.570   &   0.540   &   1.172   &   0.540      \\
200 &   0.559   &   1.046   &   0.511    &   0.542   &   1.026   &   0.507      \\
201 &   0.550   &   1.221   &   0.550    &   0.528   &   1.119   &   0.528      \\
500  &   0.538  &     1.030 &      0.507&       0.527&       1.017 &      0.504 \\
501  &   0.532   &    1.135 &      0.532&       0.518&       1.074&       0.518 \\
1000 &       0.527 & 1.022&       0.505 &       0.519&       1.012 &      0.503\\
1001 &       0.522 & 1.094 &      0.522 &     0.513&       1.052&       0.513 \\
\hline
\end{tabular}
\end{small}
\caption{Weight of the left tail $w_L=P(x\leq A)$ and the ratio of the weights of two tails $w_L/w_R$ for
 $Binom(n,p)$ distribution.
$w_L^m$ stands for the modified weight $w_L^m=P(X\leq A)/(1+P(A))$.} \label{tab:binom}
\end{table}
\end{center}

 \subsection {Hypergeometric distribution}
Consider a crosstabulation of two binary variables $A$ and $B$.
We shall refer to numbers of observations in the cell $(i,j)$ and respective probabilities as $n_{ij}$ and $p_{ij}$,
 $i,j=1,2$. The value
$n_{11}$ is the statistic of Fisher's exact test used to test for association of $A$ and $B$ given fixed margins
 $n_{1+},n_{+1},n$.  The value $n_{11}$  defines all the other entries
in a table with given margins. A parameter of primary importance is the odds ratio $\rho=p_{11}p_{22}/p_{12}p_{21}$
estimated by $\hat\rho=n_{11}n_{22}/n_{12}n_{21}$. The case
of no association $p_{ij}=p_{i+}p_{+j}$ is equivalent to $\rho=1$. Denote the expected values
$m_{ij}=\e(n_{ij})=n_{i+}n_{+j}/n$, with $E=m_{11}=\e(n_{11})$.
The number $n_{11}>E$ iff $\hat\rho>1$ .
\citet{Fisher} derived the distribution of $n_{11}$ as
$$f(n_{11}; n_{1+},n_{+1};  \rho) = \frac {{n_{1+}\choose n_{11}}{n-n_{1+}\choose n_{+1}-n_{11}}\rho^{n_{11}}}
{\sum_u {n_{1+}\choose u}{n-n_{1+}\choose n_{+1}-u}\rho^u}.$$
The null distribution (standard hypergeometric) is for $\rho = 1$. For testing $H_0 : \rho = 1$ vs $H_1 : \rho > 1$ the p-value is $p_+=\sum_{u\geq n_{11}}f(u; n_{1+},n_{+1};  \rho)$.
For $H_1 : \rho < 1$ the p-value is $p_-=\sum_{u\leq n_{11}}f(u; n_{1+},n_{+1};  \rho)$. For a two-sided test, $P_{prob}$ seems to be the p-value of choice, implemented both in R and in StatExact.\\

Sometimes other one-sided test statistics are used to test for association; they may be based on the differences of proportions in rows
or columns (e.g. $n_{11}/n_{+1}-
n_{12}/n_{+2}$) or on the $\log \hat\rho$. Nevertheless, all other possible 1-sided  tests are equivalent to
Fisher's exact test since their statistics are strictly increasing functions of $n_{11}$, as shown by \citet{Davis}.
The Fisher's exact test is also the UMPU test if the randomization is allowed \citep{Tocher}.\\

For $Hyper(x;n_{1+},n_{+1},n)$ distribution,   the full range of values $x$ for fixed margins $(n_{1+},n_{+1},n)$
is $\{ x=m_-,\cdots,m_+\}$, where $m_-=\Max(0,n_{1+}+n_{+1}-n)$ and $m_+= \Min(n_{1+},n_{+1})$.
The mode is $M=\lfloor (n_{1+}+1)(n_{+1}+1)/(n+2)\rfloor=\\
\lfloor \frac{n}{n+2}(p_{1+}(1-p_{+1})+p_{+1}(1-p_{1+})+1/n)+E\rfloor$. Therefore $\lfloor E\rfloor\leq M\leq \lfloor E+1/2\rfloor$.
When $M$ is an integer, $M-1$ and $M$ are both modes
and the mean $E\in(M-1,M)$ is unattainable. When $E$ is an integer, $M=E$.
In all cases the distance $|M-E|<1$.\\

Exact 2-sided tests for association are used when both
positive and negative associations are of interest. However, there
is ongoing controversy about how 2-sided p-values should be
constructed for the hypergeometric distribution \citep{Yates, Agresti-1977, meulepas, Dunne}.
\\

\citet{Davis} compares the p-values associated with the following 6 statistics:
$T_1=-P(n_{11}),\; T_2=|n_{11}/n_{+1}-n_{12}/n_{+2}|=N(n_{+1}n_{+2})^{-1}|n_{11}-m_{11}|,\;\\
T_3=|n_{11}/n_{1+}-n_{21}/n_{2+}|=N(n_{1+}n_{2+})^{-1}|n_{11}-m_{11}|,\;T_4=|\log(\hat\rho)|,\;\\
T_5=\sum_{ij}(n_{ij}-m_{ij})^2/m_{ij}=n^3(n_{11}-m_{11})^2(n_{1+}n_{2+}n_{+1}n_{+2})^{-1},\;\\ T_6=2\sum_{ij}n_{ij}\log(n_{ij}/m_{ij}).$
Statistic $T_1$ orders the the tables according to their probability, and corresponds to a test based on $P_{prob}$,
$T_2$ and $T_3$ are standard large-sample tests for homogeneity of proportions, $T_4$ \citep{HP-1965} rejects
for small and large values of observed log-odds ratio, $T_5$ is the Pearson's chi-square test
statistic, and $T_6$ is the likelihood ratio statistic \citep{Agresti-1977}.
It can be seen that $T_2,\;T_3$ and $T_5$ are strictly increasing functions of $|n_{11}-m_{11}|$, and therefore
the p-values for them do not differ. Further, all of the statistics $T_j, j=1,\cdots,6$ are strictly  decreasing functions of $n_{11}$ for
$n_{11}\leq m_{11}$, and strictly increasing functions of $n_{11}$ for
$n_{11}\geq m_{11}$.  \citet{Davis} further shows that the 2-sided tests $T_1,\;T_4,T_5$ and $T_6$
are not equivalent due to differing ordering of the tables at the opposite tail.

Consider the table with margins $(n_{1+}n_{2+}n_{+1}n_{+2})=(9,21,5,25)$ used as an example in \citet{Davis}.
The possible $n_{11}$ values are 0 through 5, $\e(n_{11})=1.5$, so the left tail has two tables only, for $n_{11}=0$ and 1,
with the total probability of $w_L=.521$. Tables with $n_{11}=2,\cdots,  5$ are on the right tail, the total
probability is $w_R=.479$.
 Two tails are rather close in probability.
 \citet{Davis} looks at the orderings of tables according to the increasing values of test statistics, as follows:
 $$\begin{array}{ll}
 T_1:&1\,\;2\,\;0\,\;3\,\;4\,\;5\\
 T_4:&2\,\;1\,\;3\,\;4\,\;0\,\;5\\
 T_5:&\{1\,\;2\}\,\;\{0\,\;3\}\,\;4\,\;5\\
 T_6:&2\,\;1\,\;3\,\;0\,\;4\,\;5
 \end{array}$$

 Due to monotonicity of all statistics $T_j,j=1,\cdots,6$ at both sides of the mean $m_{11}$, the
 conditional p-values for all 6 statistics do not differ (Lemma 2). Therefore all 6 2-sided tests are equivalent. This is the main advantage
 of the conditional p-value for hypergeometric distribution. Fisher's exact test is usually superseded by the chi-square test for large cell numbers.
 Equivalence of these two tests is of practical importance, for example when testing for linkage disequilibrium in genetics.\\

The probabilities of the 6 tables along with their one-sided p-values, $P_{prob}$ and $P_C$ values are
 given in columns 2-5 of Table 2.
 Conditional p-values are very close to doubled 1-sided p-values.
The second set of tables in Table 2
 corresponds to margins $(n_{1+}n_{2+}n_{+1}n_{+2})=(9,31,5,35)$ .
Here the left tail probability is 0.689, and the thin right tail has probability 0.311. The probabilities
and the p-values are given in columns 6-9. Here the inflation of the conditional p-values on the right tail is more prominent.

 \begin{table}[h]
   \label{tab:davis}
   \begin{small}
   \begin{tabular}{|c||r|r|r|r||r|r|r|r|}
     \hline
    $n_{11}$   & $P(n_{11})$ & $p_{1-sided}$& $P_{prob} $ & $P_C$ & $P(n_{11})$ & $p_{1-sided}$& $P_{prob} $ & $P_C$\\
    \hline
    \hline
 0 & .143  & .143 & .286 & .274&.258&.258&.570&.374\\

 1 & .378  & .521  & 1 & 1 &.430&.689&1&1 \\

 2 & .336  & .479 &  .622 & 1&.246&.311&.311&1   \\

 3 & .124  & .143  & .143   & .299 &.059&.065&.065&.209 \\

 4 & .019  &  .019 & .019    & .040&.006&.006&.006&.028  \\
 5 & .001  &  .001 & .001   & .002 &.0002&.0002&.0002&.0006\\
   \hline
     \hline
     \end{tabular}
   \end{small}
\caption{\small
6 possible tables, their probabilities and various p-values
 for Fisher's exact test for a table with margins $(n_{1+}n_{2+}n_{+1}n_{+2})=(9,21,5,25)$
  are given in columns 2-5. The same information for a table with
margins $(n_{1+}n_{2+}n_{+1}n_{+2})=(9,31,5,35)$ is given in columns 6-9. }
\end{table}

\section{Discussion}
Two-sided testing in non-symmetric distributions is not straightforward. The UMPU tests are not implemented in the mainstream software packages even for continuous problems, and require randomization in the discrete case. The non-asymptotic GLR tests are also not implemented, and are, in general, biased, \citet{BarLev2002}. At the same time the two-sided tests are the staple in all applications. An importance of a conceptually and computationally simple approach to two-sided testing is self-evident.

The conditional 2-sided p-value $P_C$ introduced in Section 1 is closely related
to doubled p-value and has an intuitive appeal. Its use is advocated
for both continuous and discrete distributions. An important  advantage of this p-value
is that equivalent 1-sided tests are transformed into  $P_C$-equivalent
  2-sided tests. This helps to resolve the ongoing controversy about which 2-sided tests should be used for the association
  in 2 by 2 tables.\\

   The properties of this p-value compare favorably to the doubled p-value and to the minimum likelihood p-value $P_{prob}$ , the main two implemented options in  statistical tests for non-symmetric distributions. For the variance test,
  the bias of the $P_C$-based test is smaller than the bias of the standard equal tails test based on the doubled p-value, and much smaller than the bias of the $P_{prob}$-based test.  For the considerably unbalanced sample sizes, the $P_C$-based test is also less biased than the equal tails  F-test of the equality of variances.

    We did not compare the power and the bias of the resulting tests for the binomial and the hypergeometric cases. This is
  difficult to do for tests at different levels without recourse to randomisation. For asymptotically normal
  tests, both p-values should result in asymptotically UMPU tests, 
  though the minimum likelihood p-value may require more stringent conditions to ensure the convergence
  of the density to normal density. The proof of these statements is a matter for further research.

   Another open question is
  which version $P_C^A(x)$ or $P_C^{A(m)}(x)$ should be used for an attainable value of $A$. Motivation for  $P_C^{A(m)}(x)$ is less clear, it also results in a more conservative test on top of the inescapable conservativeness due to discrete distribution.

 \citet{GP-1975} consider a large number of 2-sided p-values and find them lacking. They recommend reporting one-tailed p-value with the direction of the observed departure
  from the null hypothesis. In this spirit, the conditional p-value conditions on this direction.



\bibliographystyle{apalike}
\bibliography{twisided}

\end{document}